\magnification 1250
\mathcode`A="7041 \mathcode`B="7042 \mathcode`C="7043
\mathcode`D="7044 \mathcode`E="7045 \mathcode`F="7046
\mathcode`G="7047 \mathcode`H="7048 \mathcode`I="7049
\mathcode`J="704A \mathcode`K="704B \mathcode`L="704C
\mathcode`M="704D \mathcode`N="704E \mathcode`O="704F
\mathcode`P="7050 \mathcode`Q="7051 \mathcode`R="7052
\mathcode`S="7053 \mathcode`T="7054 \mathcode`U="7055
\mathcode`V="7056 \mathcode`W="7057 \mathcode`X="7058
\mathcode`Y="7059 \mathcode`Z="705A
\def\spacedmath#1{\def\packedmath##1${\bgroup\mathsurround
=0pt##1\egroup$} \mathsurround#1
\everymath={\packedmath}\everydisplay={\mathsurround=0pt}} 
\def\nospacedmath{\mathsurround=0pt
\everymath={}\everydisplay={} } \spacedmath{2pt}
\def\up#1{\raise 1ex\hbox{\smallf@nt#1}}
\def\tx{\kern-1.5pt -}

\def\cqfd{\kern 2truemm\unskip\penalty 500\vrule height
4pt depth 0pt width 4pt\medbreak} \def\carre{\vrule height
4pt depth 0pt width 4pt} \def\virg{\raise .4ex\hbox{,}}  
\def\no{n\up{o}\kern 2pt}
\def\ind{\par\hskip 1cm\relax}

\def\pc#1{\tenrm#1\sevenrm}
\def\Res{\mathop{\rm Res}\nolimits}
\def\Card{\mathop{\rm Card}\nolimits}
\def\Div{\mathop{\rm Div}\nolimits}
\font\eightrm=cmr8
\input amssym.def
\input amssym
\catcode`\@=11
\mathchardef\dabar@"0\msafam@39
\def\lda{\mathrel{\dabar@\dabar@\dabar@\dabar@\dabar@\dabar@
\mathchar"0\msafam@4B}}
\catcode`\@=12
\def\dfl#1{\buildrel {#1}\over {\lda}}
\def\diagramme#1{\def\normalbaselines{\baselineskip=0truept
\lineskip=4truept\lineskiplimit=1truept}   \matrix{#1}}
\def\vfl#1#2{\llap{$\scriptstyle#1$}\left\downarrow\vbox to
6truemm{}\right.\rlap{$\scriptstyle#2$}}

\parindent=0cm
\vsize = 25truecm
\hsize = 16truecm
\voffset = -.5truecm
\baselineskip15pt
\overfullrule=0pt\frenchspacing
\centerline{\bf Endomorphisms of hypersurfaces and other manifolds}
\smallskip
\smallskip \centerline{Arnaud {\pc BEAUVILLE}} 
\vskip1.2cm

\ind We prove in this note the following result:\smallskip 
{\bf Theorem}$.-$ {\it A smooth hypersurface $X$ of dimension
$n\ge 2$ and degree $d\ge 3$ admits no endomorphism of degree}
$>1$.\par
\ind  Since the case
of quadrics is treated in [PS], this settles the question of
endomorphisms of hypersurfaces. We prove the theorem in Section 1, 
using a simple but efficient trick devised by Amerik, Rovinsky and
Van de Ven [ARV]. In Section 2 we collect some general results on
endomorphisms of projective manifolds; we prove in particular that
{\it ramified} endomorphisms occur only on varieties of 
Kodaira dimension $-\infty$. This leads naturally to ask the
existence problem for Fano manifolds; we will settle this question
for surfaces.
\vskip0.4truecm
\ind  {\eightrm  I am indebted to I. Dolgachev for bringing the
problem to my attention.}

\vskip0.8cm
{\bf 1. Hypersurfaces}
\ind The proof of the theorem is based on the following result,
which appears essentially in
[ARV]: \smallskip  {\bf Proposition 1}$.-$ {\it Let $X$ be a
submanifold of ${\bf P}^N$, of dimension $n$, and let
$f:X\rightarrow X$ be an endomorphism of $X$ such that
$f^*{\cal O}_X(1)={\cal O}_X(m)$ for some integer $m\ge 2$.
Then}$$c_n(\Omega^1_X(2)\le 2^n\,\deg(X)\ .$$
\ind Let us sketch the proof following [ARV]. We first observe
that the sheaf $\Omega^1_{{\bf P}^N}(2)$ is spanned by its global
sections; therefore $\Omega^1_X(2)$, which is a quotient of
$\Omega^1_{{\bf P}^N}(2)^{}_{\,|X}$, is also spanned by its global
sections. Let $\sigma $ be a general section of
$\Omega^1_X(2)$; then $\sigma $ and its pull-back
$f^*\sigma \in H^0(X,\Omega^1_X(2m))$ have isolated zeroes
[\hbox{ARV}, lemma 1.1]. Counting these zeroes gives
$$c_n(\Omega^1_X(2m))\ge \deg(f)\ c_n(\Omega^1_X(2)) \ .$$
Since  $\deg(f)=m^n$ we get $c_n(\Omega^1_X(2))\le
m^{-n}\,c_n(\Omega^1_X(2m))$. Replacing $f$ by $f^k$ we obtain
this inequality for $m$ arbitrarily large; therefore 
$$c_n(\Omega^1_X(2))\le \lim_{m\rightarrow
\infty}m^{-n}\,c_n(\Omega^1_X(2m))=2^n\,\deg(X)\ .\ \ \carre$$
\bigskip
{\it Proof of the Theorem} :  We first discuss the case $n\ge 3$. 
Then $b_2(X)=1$, so that the condition $f^*{\cal O}_X(1)={\cal
O}_X(m)$ is automatic. In view of the Proposition we just have to
prove that $c_n(\Omega^1_X(2)> 2^n\,\deg(X)$. From the exact
sequences
$$\nospacedmath\displaylines{0\rightarrow \Omega^1_{{\bf
P}^{n+1}}(2)^{}_{|X}\longrightarrow {\cal
O}^{}_X(1)^{n+2}\longrightarrow {\cal O}_X(2)\rightarrow 0\cr
0\rightarrow {\cal O}_X(2-d)\longrightarrow \Omega^1_{{\bf
P}^{n+1}}(2)^{}_{|X}
\longrightarrow \Omega^1_X(2)
\rightarrow 0}$$we get $c(\Omega^1_X(2))=(1+h)^{n+2}(1+2h)^{-1}
(1+(2-d)h)^{-1} $, so that
$$c_n(\Omega^1_X(2))=d\,\Res_{0} \omega\quad{\rm with}\quad
\omega={(1+x)^{n+2}\over x^{n+1}(1+2x) (1+(2-d)x)}\ \cdot $$
Straightforward computations give
$$\Res_{\infty}\omega={1\over 2(d-2)}\qquad \Res_{-{1\over
2}}\omega={(-1)^{n+1}\over 2d}\qquad \Res_{{1\over d-2}}\omega=
{-(d-1)^{n+2}\over d(d-2)}\ \virg$$  hence, by the residue theorem,
$$c_n(\Omega^1_X(2))={2(d-1)^{n+2}-d+(-1)^{n}(d-2)\over 2(d-2)}\ 
\cdot $$ Using $(d-1)^2=d(d-2)+1$ we get $c_n(\Omega^1_X(2))>
d(d-1)^n\ge  d\,2^n$, hence the result in this case.
\smallskip\ind For the  case $n= 2$, we observe that the result
is straightforward when $K_X$ is ample or trivial (see
Proposition 2 below); therefore it only  remains to prove it
for  cubic surfaces. This can be easily done with the above
method, but we will deduce it from 
 the more general case of Del
Pezzo surfaces (Proposition 3).\cqfd
\bigskip
{\it Remark}$.-$ The same method applies (with some work) to
complete intersections of multidegree $(d_1,\ldots,d_p)$ in
${\bf P}^{n+p}$, {\it provided one of the $d_i$ is} $\ge 3$.
On the other hand it does not work in general for complete
intersection of quadrics.
\vskip.8cm
{\bf 2. Other manifolds}\smallskip \ind Let $X$ be a compact 
manifold, and let  $f$ be an endomorphism of $X$ degree $>1$; by
this we mean that $f$ is generically finite (or
equivalently surjective). If
$X$ is projective (or more generally K\"ahler), $f$ is actually
finite: otherwise it contracts  some curve $C$ to a point, so that
the class of $[C]$ in
$H^*(X,{\bf Q})$ is mapped to $0$ by $f_*$. This contradicts the
following remark:\smallskip {\it Lemma$\,1\,.-$ Let $d=\deg f$. The
endomorphisms
$f^*$ and $d^{-1} f_*$ of $H^*(X,{\bf Q})$ are inverse of each
other}.\par
\ind This follows from the formula $f_*f^*=d\,{\rm
Id}$.\cqfd\medskip

\ind The existence of an  endomorphism of degree $>1$ has strong
implications on the Kodaira dimension  of $X$:\smallskip
{\bf Proposition 2}$.-$ {\it Let
$X$ be a compact manifold, with an endomorphism $f$ of degree
$>1$}.\ind a) {\it The Kodaira dimension}
$\kappa(X)$ {\it is} $<\dim(X)$.\ind b) {\it If
$\kappa(X)\ge 0$, $f$ is \'etale}. \par {\it Proof} : a)
follows for instance from [KO]; let us give the proof for
completeness. Consider the pluricanonical maps
$\varphi_m: X\dasharrow {\bf P}(H^0(X,mK_X)$ associated to the
linear systems $|mK_X|\ (m\ge 1)$.  
 The pull-back map
$f^*:H^0(X,mK_X)\rightarrow
H^0(X,mK_X)$ is injective, and
therefore bijective; we have a commutative diagram:
$$\diagramme{X & \dfl{\varphi_m} & {\bf
P}(H^0(X,mK_X))&\cr
\vfl{f}{} &&\vfl{\wr\kern-3pt}{\kern-2pt{\bf P}(f^*)}&\cr
X & \dfl{\varphi_m} & {\bf
P}(H^0(X,mK_X))&.}$$
In particular, we see that $f$ {\it induces an automorphism of}
$\varphi_m(X)$. If $\dim \varphi_m(X)=$ $\dim X$ this implies
$\deg f=1$.
\ind b) Let $m$ be a
positive integer such  that the linear system
$|mK_X|$ is non-empty. Let $F$ be the
fixed divisor of this system, and $|M|$ its moving part, so that 
$mK_X\equiv F+M$. The Hurwitz formula reads $K_X\equiv
f^*K_X+R$, where $R$ is the ramification divisor of $f$; this gives
$$F+M\equiv (f^*F+mR)+f^*M\ .$$
\ind In particular we have $h^0(f^*M)\le h^0(M)=h^0(mK_X)$; since
the pull-back map $f^*:H^0(X,M)\rightarrow H^0(X,f^*M)$ is
injective, we  get $h^0(f^*M)= h^0(mK_X)$, which means that
$|f^*M|$ is the moving part of $|mK_X|$ and $f^*F+mR$ its fixed
part. Thus
$$F=f^*F+mR$$in the divisor group $\Div(X)$ of $X$. Let
$\nu:\Div(X)\rightarrow {\bf Z}$ be the homomorphism which takes
the value 1 on each irreducible divisor. 
Since $\nu(f^*F)\ge\nu(F)$, the above equality is possible only if
$R=0$: we conclude that $f$ is \'etale.\cqfd
\medskip
\ind Every Kodaira dimension $<\dim X$ can indeed occur, as
shown by the varieties  $V\times A$, where $A$ is an abelian
variety. It seems possible that all examples with
$\kappa(X)\ge 0$ are of this type, up to an \'etale covering
and perhaps some birational transformation. We can make this
precise for surfaces:
\smallskip 
{\bf Proposition 3}$.-$ {\it Let $S$ be a projective surface
with $\kappa(S)\ge 0$, admitting an endomorphism  of degree
$>1$. Then $S$ is  an abelian surface or a quotient
  $(E\times C)/G$, where $E$ is an elliptic curve, $C$ a curve
of genus $\ge 1$, and $G$  a finite group of automorphisms of
$E$ and $C$ acting freely on $E\times C$.}\par
{\it Proof} : We first observe that the surface $S$ is
minimal: if $E$ was an exceptional curve on $S$, its pull back
$f^{-1} (E)$ would be a disjoint union of exceptional curves
$E_1,\ldots,E_d$ on $S$, with $d=\deg f$. These curves would have
different classes in $H^2(S,{\bf Q})$ mapping to the same class
$[E]$ under
$f_*$, contrary to Lemma 1.
\ind By Proposition 2,  $f$ is \'etale; this implies that the
topological Euler number $e(S)$ is zero. Also we have
$\kappa(S)=0$ or $1$. In the first case, the classification
of surfaces shows that $S$ is abelian or bielliptic -- that
is, of the form  $(E\times C)/G$ with both $E$ and $C$
elliptic.
If $\kappa(S)=1$, $S$ admits an elliptic fibration
$S\rightarrow B$; since $e(S)=0$ the fibres of $f$ are (possibly
multiple) smooth elliptic curves. It is then well-known that
$S$ is isomorphic to a quotient $(E\times C)/G$ (see e.g. [B],
chap. VI).\cqfd\medskip
\ind Let us now turn  to {\it ramified} endomorphisms. By
Proposition 2 we must consider manifolds with
$\kappa(X)=-\infty$; a natural place to look at is 
Fano manifolds. For surfaces we have a complete
answer:\smallskip 
{\bf Proposition 4}$.-$ {\it A Del Pezzo surface $S$
admits an  endomorphism of degree $>1$ if and only if $K_S^2\ge
6$}.\par {\it Proof} : a) A Del Pezzo surface of degree $\ge
6$ is isomorphic to ${\bf P}^1\times {\bf P}^1$ or ${\bf P}^2$
blown up at some of the points $(1,0,0),\ (0,1,0),\ (0,0,1)
$.  The first case is trivial; in the second case,  
the endomorphisms $(X,Y,Z)\mapsto (X^p,Y^p,Z^p)$ of ${\bf
P}^2$ extend to the blown up surface.
\ind b) Let us now consider a Del Pezzo surface $S$ with an 
endomorphism $f:S\rightarrow S$ of degree $d>1$. Let $E$ be an
exceptional curve on $S$, $F=f(E)$, and $\delta $ the degree of
$f_{|E}:E\rightarrow F$. We have $f_*E=\delta F$ and therefore 
$f^*F\equiv {d\over \delta }\,E$ (Lemma 1). Taking squares gives
$F^2=-{d\over \delta^2 }$. Because of the genus formula
$C^2+C.K=2g(C)-2$,  the only curves with
negative square on a Del Pezzo surface are the exceptional
ones. Thus $F$ is exceptional,
$d=\delta ^2$ and $f^*F\equiv \delta E$; since the right hand
side does not move, this is an equality of divisors. It means
that $f$ is ramified along $E$ with rami\-fi\-cation index $\delta
$. In other words, if we denote by ${\cal E}$ the
(finite) set of exceptional curves on $S$ and by $R$  the
ramification divisor  of $f$, we have
$\displaystyle R=\sum_{E\in{\cal E}}(\delta -1)E\ +Z$, where
$Z$ is an effective divisor. Intersecting  with
$-K_S$ gives
$$-K_S\cdot R\ \ge\ (\delta -1)\Card ({\cal E})\ .$$
\ind For each $E\in{\cal E}$ we have $f^*K_S\cdot E=K_S\cdot
f_*E=\delta\, K_S\cdot F=-\delta $, and therefore
$(f^*K_S-\delta\, K_S)\cdot E=0$. We can assume that ${\cal
E}$ spans the Picard group of $S$ (this holds as soon
as $K_S^2\le 7$), thus $f^*K_S\equiv \delta\, K_S$. Then the
Hurwitz formula $K_S\equiv f^*K_S+R$  gives $R\equiv (\delta
-1)(-K_S)$, so that the above inequality becomes $K_S^2\ge
\Card({\cal E})$. This is impossible for $K_S^2\le 5$, as the
surface $S$ contains then at least 10
exceptional curves.\cqfd\medskip
\ind For Fano threefolds we know the answer in the case
$b_2=1$, as a consequence of the more general results of [A]
and [ARV]: the only Fano threefold with $b_2=1$ admitting an
endomorphism of degree $>1$ is ${\bf P}^3$. Their methods
apply to some other Fano threefolds, but the general case
seems to  require new techniques.
\vskip2cm
\centerline{REFERENCES} \vglue15pt\baselineskip13pt
\def\num#1{\smallskip\item{\hbox to\parindent{\enskip [#1]\hfill}}}
\parindent=1.1cm 
\num{A} E. {\pc AMERIK}: {\sl Maps onto certain Fano
threefolds}. Doc. Math. {\bf 2} (1997), 195--211.
\num{ARV} E. {\pc AMERIK}, M. {\pc ROVINSKY}, {\pc A. VAN DE} {\pc
VEN}: {\sl A boundedness theorem for morphisms between threefolds}.
Ann. Inst. Fourier {\bf 49} (1999), 405--415.
\num{B} A. {\pc BEAUVILLE}: {\sl 	Surfaces alg\'ebriques complexes.} 
Ast\'erisque {\bf 54}  (1978).
\num{KO} S. {\pc KOBAYASHI}, T. {\pc OCHIAI}: {\sl Meromorphic
mappings onto compact complex spaces of general type}. Invent.
math. {\bf 31} (1975), 7--16.
\num{PS} K. {\pc PARANJAPE}, V. {\pc SRINIVAS}: {\sl Self maps of
homogeneous spaces}. Invent. math. {\bf 98} (1989), 425--444.
\vskip1cm\baselineskip11pt
\font\eightrm=cmr8
\font\sixrm=cmr6
\def\pc#1{\eightrm#1\sixrm}
\hfill\vtop{\eightrm\hbox to 5cm{\hfill Arnaud {\pc BEAUVILLE}\hfill}
\hbox to 5cm{\hfill Laboratoire J.-A. Dieudonn\'e\hfill}
\hbox to 5cm{\sixrm\hfill UMR 6621 du CNRS\hfill}
\hbox to 5cm{\hfill {\pc UNIVERSIT\'E DE}  {\pc NICE}\hfill}
\hbox to 5cm{\hfill  Parc Valrose\hfill}
\hbox to 5cm{\hfill F-06108 {\pc NICE} Cedex 02\hfill}}

\end